\documentclass[]{article}
\input{amssym}
\input{epsf}
\begin{document}
\def\be{\begin{eqnarray*}}
\def\ee{\end{eqnarray*}}
\def\di{\displaystyle}
\title{{\bf Lie symmetry analysis of the Grad-Shafranov equation}}
\author{ Mehdi Nadjafikhah\thanks{School of
Mathematics, Iran University of Science and Technology, Narmak-16,
Tehran, I.R.Iran. e-mail:~m\_nadjafikhah@iust.ac.ir}\and Parastoo
Kabi-Nejad\thanks{e-mail:~parastoo\_kabinejad@iust.ac.ir}}
\date{ }
\maketitle
\begin{abstract}
The theory of plasma physics offers a number of nontrivial
examples of partial differential equations, which can be
successfully treated with symmetry methods. We propose the
Grad-Shafranov equation which may illustrate the reciprocal
advantage of this interaction between plasma physics and symmetry
techniques. A symmetry classification of the Grad-Shafranov
equation with two arbitrary functions $F(u)$ and $G(u)$ of the
unknown variable $u=u(x,t)$ is given. The optimal system of
one-dimensional subalgebras is performed. This latter provides a
process for building new solutions for the equation.
\end{abstract}
\medskip \noindent {\bf A.M.S. 2000 Subject Classification:} 34C14, 35J05, 70G65.

\section{Introdutction}
The equation we are going to examine contains indeed two
arbitrary functions $F(u), G(u)$ of unknown variable $u=u(x,t)$,
and the goal is now to perform the symmetry classification of
this equation. i.e. to find those $F, G$ for which the equation
admits nontrivial symmetries. In general, the symmetry
propoerties of an equation may strongly depend on the choice of
the arbitrary functions involved.

The PDE we want to consider is
\begin{eqnarray}
u_{xx}+a\frac{1}{x}u_x+u_{tt}=x^pF(u)+G(u)
\end{eqnarray}
with $a=-1, b=1$ this equation is known in plasma physics as the
Grad-Shafranov equation( see \cite{[56]}) and describes the
magnetohydrodynamic force balance in a  confined toroidal plasma.
In this context, $u$ is the so-called magnetic flux variable, $x$
is a radial variable, then $x\geq 0$, while the two arbitrary flux
functions $F(u), G(u)$ are related to the plasma pressure and
current density profiles.

\medskip\noindent This equation admits well known properties and a rich
literature is devoted to it. After some comments on the general
peculiarities of the equation and of its symmetry properties, we
shall provide the algebra of its exact Lie point symmetries and
probe some properties of its algebra. Next, we shall study
one-dimensional optimal system of  subalgebras.

\medskip\noindent According to the standard definitions and procedure (see
\cite{[1]}), we can now look for the equivalence group;
preliminary, we look for the full groups of the equation
(\ref{eq:2}), i.e. the intersection of all groups admitted by
(\ref{eq:2}) for any arbitrary choice of $F$ and $G$. It turns
out that this kernel is almost trivial: it contains indeed only
the transformation of the variable t.

\medskip\noindent First of all, we exclude from our classification the case
$a=0$ and $p=0$ (or $a=0$ and $F(u)\equiv0$) because in this case
our equation coincides with well known nonlinear Laplace equation.

\medskip\noindent Now, we investigate the case, $F=1$ and $G=0$ for
the Grad-Shafranov equation.
\section{Symmetry methods}
Let a partial differential equation consists $p$ independent and
$q$ dependent variables. The one-parameter Lie group of
transformations
\begin{eqnarray}
\bar{x_i}=x_i+s\xi^i(x,u)+O(s^2);\qquad
\bar{u_\alpha}=u_\alpha+s\varphi^\alpha(x,u)+O(s^2), \label{eq:2}
\end{eqnarray}
where $\xi^i=\frac{\partial\bar{x_i}}{\partial
s}|_{s=0},i=1,\cdots, p,$ and
$\varphi^\alpha=\frac{\partial\bar{u_\alpha}}{\partial s}|_{s=0},
\alpha=1,\cdots, q,$ are given. The action of the Lie group can
be recovered from that of its associated infinitesimal generators.
We consider general vector field
\begin{eqnarray}
Y=\sum_{i=1}^p\xi^i(x,u)\frac{\partial}{\partial x_i}+
\sum_{\alpha=1}^q\varphi^\alpha(x,u)\frac{\partial}{\partial
u^\alpha} \label{eq:3}
\end{eqnarray}
on the space of independent and dependent variables. Therefore,
the characteristic of the vector field $Y$ given by (\ref{eq:3})
is the function
\begin{eqnarray}
Q^\alpha(x,u^{(1)})=\varphi^\alpha(x,u)-\sum_{i=1}^p\xi^i(x,u)\frac{\partial
u^\alpha}{\partial x_i},\quad \alpha=1,\cdots,q. \label{eq:4}
\end{eqnarray}
The second prolongation of the infinitesimal generator
\begin{eqnarray}
X=\xi^1(x,t,u)\frac{\partial}{\partial
x}+\xi^2(x,t,u)\frac{\partial}{\partial
t}+\varphi(x,t,u)\frac{\partial}{\partial u} \label{eq:5}
\end{eqnarray}
is the following vector field
\begin{eqnarray}
X^{(2)}&=&X+\varphi^x\frac{\partial}{\partial
u_x}+\varphi^t\frac{\partial}{\partial u_t}
+\varphi^{xx}\frac{\partial}{\partial
u_{xx}}+\varphi^{xt}\frac{\partial}{\partial u_{xt}}+
\varphi^{tt}\frac{\partial}{\partial u_{tt}}\nonumber\label{eq:6}
\end{eqnarray}
with coefficients
\begin{eqnarray}
\varphi^x&=&D_x Q+\xi^1u_{xx}+\xi^2u_{xt},\\
\varphi^t&=&D_t Q+\xi^1u_{xt}+\xi^2u_{tt},\nonumber\\
\varphi^{xx}&=&D^2_x Q+\xi^1u_{xxx}+\xi^2u_{xxt},\nonumber\\
\varphi^{xt}&=&D_xD_t Q+\xi^1u_{xxt}+\xi^2u_{xtt},\nonumber\\
\varphi^{tt}&=&D^2_t Q+\xi^1u_{xtt}+\xi^2u_{ttt},\nonumber
\label{eq:7}
\end{eqnarray}
where the operators $D_x$ and $D_t$ denote the total derivatives
with respect to $x$ and $t:$

\medskip\noindent By theorem 6.5. in \cite{[17]}, the vector field $X$ is a
one parameter of Grad-Shafranov equation if and only if
\begin{eqnarray}
\Pr^{(2)}X(u_{xx}-\frac{1}{x}u_x+u_{tt}-x^2)=0\qquad whenever
\qquad u_{xx}-\frac{1}{x}u_x+u_{tt}-x^2=0.
\end{eqnarray}
So, we apply the criterion of infinitesimal invariance in order to
determine symmetries of the Grad-Shafranov equation\cite{[15]}.
Therefore, the infinitesimal symmetry criterion is
\begin{eqnarray}
\varphi^{xx}-\frac{1}{x^2}\xi^1u_x-\frac{1}{x}\varphi^x+\varphi^{tt}-2x\xi^1=0\label{eq:8}
\end {eqnarray}
Substituting the formulas (\ref{eq:7}) into (\ref{eq:8}), we are
left with a polynomial equation involving the various derivatives
of $u$ whose coefficients are certain derivatives of $\xi^1$,
$\xi^2$ and $\varphi$ only depend on $x, t, u$ we can equate the
individual coefficients to zero, leading to the complete set of
defining equations:
\begin{eqnarray}
\xi^2_u&=&0,
\xi^2_{tt}=-\frac{\xi^2_x}{x},\ \ \ \xi^2_{tx}=0,\ \ \ \xi^2_{xx}=\frac{\xi^2_x}{x},\\
\varphi_{tt}&=&\frac{4\xi^2_t
x^3+\varphi_x-\varphi_ux^3-\varphi_{xx}x}{x},\\
\varphi_{tu}&=&-\frac{\xi^2_x}{2x},\ \ \ \varphi_{uu}=0,\ \ \
\varphi_{ux}=0,\ \ \ \xi^1(x,t,u)=\xi^2_tx. \label{eq:9}
\end{eqnarray}
By solving this system of PDEs, we state the following theorem.
\paragraph{\bf Theorem}
{\it The Lie algebra ${\goth g}$ of the symmetry group $G$
associated to the Grad-Shafranov equation is generated by the
vector fields
\begin{eqnarray}
X_1&=&\partial_t,\nonumber\\
X_2&=&x\partial_x+t\partial_t+\frac{x^4}{2}\partial_u,\nonumber \\
X_3&=&(-\frac{x^4}{8}+u)\partial_u,\\
X_4&=&tx\partial_x+\frac{1}{2}(t^2-x^2)\partial_t+\frac{1}{16}t(7x^4+8u)\partial_u,\nonumber\\X_5&=&\psi(x,t)\partial_u.\nonumber
\label{eq:11}
\end{eqnarray}}
The commutation relations of the 4-dimensional Lie algebra ${\bf g
}$ spanned by the vector fields $X_1, X_2, X_3, X_4$ are shown in
the following table.
\begin{table}[h]
\caption{Commutation relations satisfied by infinitesimal
generators in (33).}
\begin{tabular*}{\textwidth}{@{}l*{15}{@{\extracolsep{0pt plus
12pt}}l}} \hline
$[\,,\,]$&$X_1$ &$X_2$  &$X_3$ &$X_4$ \\
\hline
  $X_1$    & 0& $X_1$ & 0 & $X_2+\frac{1}{3}X_3$ \\
  $X_2$    & $-X_1$   & 0 & 0&  $X_4$ \\
  $X_3$    & 0 & 0    & 0&0    \\
  $X_4$    & $-X_2-\frac{1}{3}X_3$ &$-X_4$&0&0\\
\hline
\end{tabular*}
\end{table}
\section{The Lie algebra of symmetries}
In this section, we determine the structure of full symmetry Lie
algebra ${\bf g}$ of Eq.  (\ref{eq:3}).
\paragraph{{\bf Theorem}}{\it The full symmetry Lie algebra {\bf g} of
Eq.(\ref{eq:3}) has the following semidirect decomposition:
\begin{eqnarray}
{\bf g}={\bf R}\ltimes {\bf g}_1
\end{eqnarray}
where ${\bf g}_1$ is a semi-simple Lie algebra.}

\medskip\noindent{\bf Proof:} The center ${\bf z}$ of ${\bf g}$ is $Span
_{\bf R}\{X_3\}$. Therefore the quotient algebra ${\bf g}_1={\bf
g}/{\bf z}$ is $Span_{\bf R}\{Y_1,Y_2,Y_3\}$, where $Y_i=X_i+{\bf
z}$ for $i=1,2,3$. The commutator table of this quotient algebra
is given in the following table:
\begin{table}[h]
\caption{Commutation relations satisfied by infinitesimal
generators in (33).}
\begin{tabular*}{\textwidth}{@{}l*{15}{@{\extracolsep{0pt plus
12pt}}l}} \hline
$[\,,\,]$&$Y_1$ &$Y_2$  &$Y_3$  \\
\hline
  $Y_1$    & 0& $Y_1$ & $Y_2$\\
  $Y_2$    & $-Y_1$   & 0 & $Y_3$ \\
  $Y_3$    & $-Y_2$ & $-Y_3$    & 0    \\
\hline
\end{tabular*}
\end{table}

The Lie algebra ${\bf g}$ is non-solvable, because
\begin{eqnarray}
{\bf g}^{(1)}&=&[{\bf g},{\bf g}]=Span_{\bf
R}\{X_1,X_2+\frac{1}{3}X_3,X_4\},\\
{\bf g}^{(2)}&=&[{\bf g}^{(1)},{\bf g}^{(1)}]={\bf g}^{(1)}.
\end{eqnarray}
Similarly, ${\bf g}_1$ is semi-simple and non-solvable, because
\begin{eqnarray}
{\bf g}_1^{(1)}=[{\bf g}_1,{\bf g}_1]=Span_{\bf
R}\{Y_1,Y_2,Y_3\}={\bf g}_1.
\end{eqnarray}
The Lie algebra ${\bf g}$ admits a Levi decomposition as the
following  semi-direct product ${\bf g}=r\ltimes s$, where
$r=Span _{R}\{X_3\}$ is the radical of ${\bf g}$ (the largest
solvable ideal contained in {\bf g}), and $s=Span_{\bf
R}\{X_1,X_2+\frac{1}{3}X_3,X_4\}$.

\medskip\noindent The ideal $r$ is a one-dimensional subalgebra of ${\bf g}$,
therefore it is isomorphic to ${\bf R}$; Thus the identity ${\bf
g}=r\ltimes s$ reduces to ${\bf g}={\bf R}\ltimes {\bf g}_1$.
\section{Optimal system of subalgebras}
Consider a system of partial differential equations $\Delta$
defined over an open subset $M\subset X\times U \simeq {\Bbb
R}^p\times {\Bbb R}^q$ of the space of independent and dependent
variables. Let $G$ be a local group of transformations acting on
$M$. Roughly, a solution $u=f(x)$ of the system is said to be
$G$-invariant if it is left unchanged by all the group
transformations in $G$. If $G$ is a symmetry group of a system of
partial differential equations $\Delta$, then, under some
additional regularity assumptions on the action of $G$, we can
find all the $G$-invariant solutions to $\Delta$ by solving a
reduced system of differential equations, denoted by $\Delta/ G$,
which will involve fewer independent variables than the original
system $\Delta$. In general, to each $s$-parameter subgroup $H$ of
the full symmetry group $G$ of a system of differential equations
in $p>s$ independent variables, there will correspond a family of
group-invariant solutions to the system. Since there are almost
always an infinite number of such subgroups, it is not usually
feasible to list all possible group-invariant solutions to the
system. We need an effective, systematic means of classifying
these solutions, leading to an "optimal system" of group-invariant
solutions from which every other such solution can be derived.
Since elements $g\in G$ not in the subgroup $H$ will transform an
$H$-invariant solution to some other group-invariant solution,
only those solution, not so related  need be listed in our optimal
system.\newline Let $G$ be a Lie group. An optimal system of
$s$-parameter subgroups is a list of conjugacy inequivalent
$s$-parameter subgroups with the property that any other subgroup
is conjugate to precisely one subgroup in the list. Similarly, a
list of $s$-parameter subalgebras forms an optimal system if every
$s$-parameter subalgebra of $\goth{g}$ is equivalent to a unique
member of the list under some element of the adjoint
representation: $\tilde{\goth{h}}= {\rm Ad}  g(\goth{h}),g \in
G$.\newline Proposition 3.7 of \cite{[15]} says that the problem
of finding an optimal system of subgroups is equivalent to that of
finding an optimal system of subalgebras. For one-dimensional
subalgebras, this classification problem is essentially the same
as the problem of classifying the orbits of the adjoint
representation, since each one-dimensional subalgebra is
determined by a nonzero vector in $\goth{g}$. This problem is
attacked by the na\"{\i}ve approach of taking a general element
$\mathbf{V}$ in $\goth{g}$ and subjecting it to various adjoint
transformations so as to "simplify" it as much as possible. Thus
we will deal with the construction of the optimal system of
subalgebras of $\goth{g}$.

The adjoint action is given by the Lie series
\begin{eqnarray}
{\rm Ad}(\exp(\varepsilon
Y_i)Y_j)=Y_j-\varepsilon[Y_i,Y_j]+\varepsilon^2[Y_i,[Y_i,Y_j]]-\cdots,
\end{eqnarray}
where $[Y_i,Y_j]$ is the commutator for the Lie algebra,
$\varepsilon$ is a parameter and $i,j=1,2,3$.

The adjoint representation of $\goth{g}$ is listed in the
following table, it consists the separate adjoint actions of each
element of $\goth{g}$ on all other elements.
\begin{table}[h]
\caption{Adjoint relations satisfied by infinitesimal generators
in (33).}
\begin{tabular*}{\textwidth}{@{}l*{15}{@{\extracolsep{0pt plus
12pt}}l}} \hline
$[\,,\,]$&$X_1$ &$X_2$  &$X_3$ &$X_4$ \\
\hline
  $X_1$&$X_1$    & $X_2+\varepsilon X_1$    & $X_2$   & $\frac{1}{2}\varepsilon^2 X_1+\varepsilon X_2+\frac{1}{2}\varepsilon X_3+X_4$\\
 $X_2$& $\exp(-\varepsilon)X_2$    & $X_2$    & $X_3$    & $\exp(\varepsilon)X_4$   \\
 $X_3$& $X_1$    & $X_2$ & $X_3$    & $X_4$    \\
 $X_4$&
 $X_1-\varepsilon X_2\frac{1}{2}\varepsilon X_3+\frac{1}{2}\varepsilon^2 X_4$&
 $X_2-\varepsilon X_1$ & $X_3$ & $X_4$\\

\hline
\end{tabular*}
\end{table}

\noindent where the (i,j)-th entry indicating ${\rm
Ad}(\exp(\varepsilon Y_i)Y_j)$.
\paragraph{\bf{Theorem}}
{\it An optimal system of one dimensional  Lie subalgebras of
Grad-Shafranov equation is provided by those generated by
\begin{eqnarray}
&1)&X_1,\quad 2)X_2, \quad  3)X_3, \quad 4)X_3-X_1, \quad 5)X_3+X_1\\
&6)&aX_2+X_3,\quad 7)aX_1+bX_2+X_4.
\end{eqnarray}
where $a,b \in {\bf R}$ are arbitrary constants.}

\medskip\noindent{\bf Proof} Let ${\bf g}$ is the symmetry group
of Eq. (\ref{eq:3}) with adjoint representation determined in
Table 2 and
\begin{eqnarray}
X=a_1X_1+a_2X_2+a_3X_3+a_4X_4
\end{eqnarray}
is a nonzero vector field of $g$. We will simplify as many of the
coefficients of $a_i, i=1,\cdots 4$; as possible through judicious
applications of adjoint maps to X.

\medskip\noindent Case 1:

\medskip\noindent Suppose first that $a_4\neq 0$. Scaling $X$ if necessary,
we can assume that $a_4=1$. Referring to table 2, if we act on
such a $X$ by $\mathrm{Ad}(\exp(-a_3X_1))$, we can make the
coefficient of $X_3$ vanish. Thus, every one-dimensional
subalgebra generated by a $X$ with $a_4\neq 0$ is equivalent to
the subalgebra spanned by $aX_1+bX_2+X_4$ where $a,b \in {\bf R}$
are arbitrary constants. No further simplifications are possible.

\medskip\noindent Case 2:

\medskip\noindent The remaining one-dimensional subalgebras are spanned by
vectors of the above form with $a_4=0$. If $a_3\neq 0$, we can
scale to make $a_3=1$. There are two subcases.

\medskip\noindent Case 2.1:

\medskip\noindent If $a_2\neq 0$, then we can  cancel the coefficient of
$X_1$ by acting on $X$ by
$\mathrm{Ad}(\exp(-\frac{a_1}{a_3}X_1))$. So that, $X$ is
equivalent to a scalar multiple $aX_2+X_3$ for some $a\in{\bf R}$.

\medskip\noindent Case 2.2:

\medskip\noindent If $a_2=0$, we can act by adjoint map generated by $X_2$ to
arrange the coefficient of $X_1$ either +1, -1 or 0. Therefore,
any one-dimensional subalgebra spanned by X with $a_4=0, a_3=1,
a_2=0$ is equivalent to one spanned by either $X_3+X_1, X_3-X_1$
or $X_3$.

\medskip\noindent Case 3:

\medskip\noindent The remaining cases, $a_3=a_4=0$, are similarly seen to be
equivalent either to $X_2 (a_2\neq 0)$ or to $X_1
(a_2=a_3=a_4=0)$.

\medskip\noindent There is not any more possible case for
studying and the proof is complete.

\medskip\noindent In continuation, we find some group invariant solutions of the
equation(\ref{eq:52}) corresponding to 1-dimensional subalgebras
generated by $X_1, X_2, X_4$.

\medskip\noindent Consider
\begin{eqnarray}
X_2=x\partial_x+t\partial_t+\frac{x^4}{2}\partial_u
\end{eqnarray}
invariants  are
\begin{eqnarray}
C_1=\frac{t}{x},\quad C_2=-\frac{1}{8}x^4+u
\end{eqnarray}
So, group invariant solution associated to these invariants has
the form
\begin{eqnarray}
u(x,t)=\frac{(\frac{x^4}{8}+C_1)\sqrt{t^2+x^2}+C_2t}{\sqrt{t^2+x^2}}
\end{eqnarray}
In the case of $X_4$, invariants are
\begin{eqnarray}
C_1=\frac{t^2+x^2}{x},\quad C_2=\frac{8u-x^4}{8\sqrt{x}}
\end{eqnarray}
therefore, we get the associated group invariant solution
\begin{eqnarray}
\frac{\frac{1}{8}(x^5t^2+x^7)\sqrt{\frac{t^2+x^2}{x}}+(2C_2t^2+C_1)x^{5/2}+C_2(x^{9/2}+\sqrt{x}t^4)}{\sqrt{\frac{t^2+x^2}{x}}x(t^2+x^2)}
\end{eqnarray}
In the case of $X=X_2+X_3$, we obtain the following invariants
\begin{eqnarray}
C_1=\frac{t}{x},\quad C_2=\frac{8u-x^4}{8x}
\end{eqnarray}
Thus corresponding invariant solution is
\begin{eqnarray}
u(x,t)=\frac{x^4}{8}+C_2\sqrt{\frac{t^2+x^2}{x^2}}x+C_1t
\end{eqnarray}
In the case of $X=X_1+X_3$, we get the invariants
\begin{eqnarray}
C_1=x,\quad C_2=\frac{-1}{8}(x^4-8u)\exp(-t)
\end{eqnarray}
therefore, we obtain the following invariant solution associated
to these invariants.
\begin{eqnarray}
u(x,t)=x\Big(\frac{1}{8}x^3+C_1\exp(t){\rm
BesselJ}(1,x)+C_2\exp(t){\rm BesselY}(1,x)\Big)
\end{eqnarray}
where ${\rm BesselJ}$ and ${\rm BesselY}$ are the ${\rm Bessel}$
functions of the first and second kinds, respectively. They
satisfy Bessel's equation.

\medskip\noindent For the vector field $X_1=\partial_t$, global invariants
are $C_1=x, \quad C_2=u$. So, the solution of reduced equation is
\begin{eqnarray}
u=x^4+4 C_1x^2+C_2
\end{eqnarray}

\section{New admitted symmetries}

In this section, we probe the remainder cases.

\medskip\noindent If $F=\exp(2u)$ and $G=\exp(u)$, there is a new admitted
symmetry which has the form
\begin{eqnarray}
X=x\partial_x+t\partial_t-2\partial_u
\end{eqnarray}
If $F=u^{1+\frac{2}{q}}$ and $G=u^{1+\frac{1}{q}}$ for all $q\neq
0 $, there is a new admitted symmetry which is
\begin{eqnarray}
X=x\partial_x+t\partial_t-2qu\partial_u
\end{eqnarray}
If $F=1, G=u$,  new symmetries are
\begin{eqnarray}
X=(x^2+u)\partial_u,\quad X=\psi(x,t)\partial_u
\end{eqnarray}
Furthermore,  the associated invariant solution for the global
invariants $C_1=x, \quad C_2=u$ is
\begin{eqnarray}
u(x,t)=-x(x-C_2{\rm BesselI(1,x)}+C_1{\rm BesselK(1,x)})
\end{eqnarray}
where ${\rm BesselI}$  and ${\rm BesselK}$ are the modified Bessel
functions of the first and second kinds, respectively  which
satisfy the modified Bessel equation.

\medskip\noindent If $F=u, G=1$, the new symmetry has the form
\begin{eqnarray}
X=(u+\psi(x,t))\partial_u
\end{eqnarray}
In addition, the corresponding invariant solution for the global
invariants $C_1=x,\quad C_2=u$ is
\begin{eqnarray}
u(x,t)=\cosh(\frac{x^2}{2})(C_1-\frac{1}{2}{\rm
Shi}(\frac{x^2}{2}))+\sinh(\frac{x^2}{2})(C_2+\frac{1}{2}{\rm
Chi}(\frac{x^2}{2}))
\end{eqnarray}
where ${\rm Shi}$ and ${\rm Chi}$ are hyperbolic Sine integral
and hyperbolic cosine integral, respectively.
\begin{eqnarray}
{\rm Shi}(x)&=&\int^x_0\frac{\sinh(t)}{t}{\rm dt}\\
{\rm Chi}(x)&=&\gamma+\ln(x)+\int^x_0\frac{\cosh(t)-1}{t}{\rm dt}
\end{eqnarray}


\end{document}